# Results from a 2016 Pilot Survey on Math Post-docs

Amy Cohen, Professor of Mathematics Emerita, Rutgers University

The survey focused on the transition from a post-doctoral position into immediately following employment. The survey went to chairs of 14 doctoral mathematics departments, most in the northeastern quadrant of the U.S. and most at public universities. The survey asked each department for data about

1. The type of the immediately next job of each recent post-doc after leaving the department
2. The type of the immediately prior job of each faculty member hired in 2 recent and 2 earlier years
3. The preparation for future careers offered by the department to its post-docs

Partial or complete responses have been received from 11 of the 14 departments surveyed. In most cases, the department filled out an on-line survey. In one case, no survey was returned, but the names of research post-docs were posted on the department's website, and survey personnel tracked them using search engines. In another case, the responding department filled out most sections of the survey and provided a list of names of research postdocs for survey personnel to track. Since few departments were able to supply data about post-docs who completed in the early 2000's, that data was not analyzed.

Since this survey was a pilot study, its results should spark discussion but should not be seen as conclusive. Anonymity was promised to departments to encourage frank responses. Related results appear in [1, 2, ] and will appear in [3,4]. Results for post-docs in the next section are consistent with those of a supplementary AMS survey [5].

## Summary of results and commentary

*<u>Next-jobs for research post-docs</u>* (data from 11 of the 14 departments)

Data were available for 162 mathematicians leaving research post-doc positions in 2012-13, 2013-14, or 2014-15. For this survey, a research post-doc position was defined to be a multi-year but non-renewable position intended to support the transition from directed thesis research to an independent research program. The definition explicitly included not only positions formally called "post-docs" but also positions serving the same purposes under other titles such as named instructorships, named assistant professorships, and certain one-year positions anticipating renewals for up to three years. The survey was sent to department chairs, but not to directors of research centers. It is unclear how many departments reported on post-docs located in centers. Of the 162 post-docs for which there was data,

- 41 (25%) obtained tenure-track positions in doctoral departments
- 36 (22%) moved to another post-doc position
- 15 ( 9%) took full-time non-tenure-track academic jobs
- 14 ( 9%) moved to tenure-track jobs in the US at bachelors or masters institutions
- 14 ( 9%) moved to non-US academic institutions (tenure status and title unknown)
- 13 (8%)  moved to business, industry, or government
- 29 (18%) were reported in the category "other/unknown"

Even assuming optimistically that two-thirds of those moving to an additional post-doc eventually find a tenure track job in a doctoral university, the pipeline from post-docs to academic research careers appears to leak. Perhaps better tracking by the post-docs' employers would give a more encouraging picture. Certainly, the movement of post-docs in applied math into business, industry, and government contributes to the nation's STEM workforce.

The flow from one post-doc to another may result from the increasing number of post-doctoral positions since the economic downturn of 2008, the decreasing number of "tenure-eligible" jobs in doctoral departments, and the increasing reliance of academic institutions on full-time non-tenure-track (FT NTT) doctoral mathematicians to provide undergraduate instruction. While some institutions do offer FT NTT faculty teaching loads consistent with productive research activity, most do not. [Cohen, "Disruptions of the Academic Math Employment Market", *AMS Notices*, pp 57-60, October 2016]

Some post-docs in this study were already in a second post-doc and moved to a third. This may reflect the attraction of geographic mobility, or the unavailability of attractive tenure-track jobs. An employment pattern of 6 years or more of post-doctoral support in not generally attractive to early career mathematicians (male as well as female) concerned about stability and work-life balance.

*(2) Prior jobs of new faculty hired from other universities*  (data from 8 of the 14 departments).
  *For starting dates in academic years 2014-15 or 2015-16*
- 25 were hired into tenure-track positions – none directly from their PhDs.
- 70 were hired into research post-doctoral positions, 61 directly from PhDs.
- 56 were hired into other full-time NTT teaching positions, 44 directly from PhDs.

  *For starting dates in academic years 2004-05 or 2005-06*
- 18 were hired into tenure-track positions – 2 directly from their PhDs.
- 21 were hired into research post-doctoral positions, 19 directly from PhDs.
- 18 were hired into other full-time NTT teaching positions, all directly from PhDs.

The increase in new hires of tenure track faculty from 18 to 25 (about 38%) may be explained by a number of factors. Attrition among senior tenured faculty may have increased due to deaths and retirements. An increasing number of undergraduates (not all of them math majors) enroll in undergraduate mathematics courses – especially in courses beyond the freshman level. These data did not distinguish between tenure track faculty hired with immediate tenure and those only eligible for eventual consideration for tenure.

The increases in numbers of research post-docs and FT NTT doctoral teaching faculty are strikingly larger – an increase from 21 to 70 research post-docs (about 233%) and from 18 to 56 FT NTT teaching faculty (about 211%). The increase in post-doc numbers are may compensate for decreasing numbers of tenure-eligible jobs. The increase in FT NTT faculty may result not only from budget considerations but also from efforts to improve graduation rates by using more effective teaching strategies in smaller classes.

*(3) Preparation offered to post-docs related to future careers* (data from 8 of the 14 departments)
- *General orientation.*
    - In addition to orientation sessions, most offered mentoring or formal supervision.
- *Strengthening mathematical knowledge for research.*
    - 7 encouraged attendance at one or more of: courses, colloquia, seminars, or conferences.
    - A few reported sponsoring lectures on applications of mathematics.
    - None reported help with finding internships.
- *Strengthening knowledge for teaching mathematics*
    - 7 reported observation of teaching, collecting student evaluations, and sharing feedback.
    - Some offered courses or seminars on topics such as use of technology, construction and grading of exams, or active student engagement in classes.
    - None tracked an individual post-doc's students' performance in later courses.
- *Strengthening communication skills*
    - 7 reported encouraging post-docs to talk at seminars, colloquia, and/or meetings.
    - Some provided guidance for writing papers and proposals.
    - 2 provided guidance for job talks and/or interviews.
    - 1 provided guidance for refereeing and reviewing.
    - None addressed communication with general audiences.

Most post-docs move to employment (academic and non-academic) in which they communicate with individuals who are not mathematicians. Large universities face pressure to increase student success, especially in entry level and STEM courses. As a recent PCAST report suggested, mathematicians do not have a good reputation for communicating effectively with students or other non-mathematical audiences.